\documentclass[10pt]{amsart}
\usepackage{amsfonts,amssymb,amscd,amsmath,enumerate,verbatim,calc}
\usepackage{mathrsfs}
\numberwithin{equation}{section}
\newtheorem{theorem}{Theorem}[section]
\newtheorem{corollary}[theorem]{Corollary}
\newtheorem{lemma}[theorem]{Lemma}
\newtheorem{proposition}[theorem]{Proposition}
\newtheorem{definition}[theorem]{Definition}

\newtheorem{example}[theorem]{Example}

 \DeclareMathOperator{\End}{End}

\begin{document}

\title[  Complex and K\"{a}hler structures on hom-Lie algebras]
 { Complex and K\"{a}hler structures on hom-Lie algebras }

\bibliographystyle{amsplain}

\author[Esmaeil Peyghan and Leila Nourmohammadifar]{ E. Peyghan and L. Nourmohammadifar}
\address{Department of Mathematics, Faculty of Science, Arak University,
Arak, 38156-8-8349, Iran.}
\email{e-peyghan@araku.ac.ir,\ l-nourmohammadifar@phd.araku.ac.ir}





\begin{abstract}
Complex and Hermitian structures on hom-Lie algebras are introduced and some examples of these structures are presented. Also, it is shown that there not exists a proper complex (Hermitian) home-Lie algebra of dimension two. Then using a hom-left symmetric algebra, a phase space is provided and then a complex structure on it, is presented. Finally, the notion of K\"{a}hler hom-Lie algebra is introduced and then using a K\"{a}hler hom-Lie algebra, a phase space is constructed. 
\end{abstract}

\maketitle







\section{Introduction}
Studying the deformations of the Witt and Virasoro algebras, lead Hartwig, Larsson and Silvestrov to introduce the notion of hom-Lie algebras \cite{HLS}. Indeed, some $q$-deformations of the
Witt and the Virasoro algebras have the structure of a hom-Lie algebra \cite{HLS, H}. Based on the close relation between the discrete, deformed vector fields and differential calculus, this algebraic structure plays important role in research fields \cite{HLS, LS}.

Differential-geometric structures play important role in the study of complex geometry. After Kodaira, K\"{a}hler structures became central in the study of deformation theory and the classification problems. Recall that an almost complex structure $J$ on a manifold $M$ is a linear complex structure (that is, a linear map which squares to -1) on each tangent space of the manifold, which varies smoothly on the manifold. A complex structure is essentially an almost complex structure with an integrability condition, and this condition yields an unitary structure ($U(n)$ structure) on the manifold. An almost Hermitian structure is a pair $(J,g)$ of an almost complex structure $J$ and a pseudo-Riemannian metric $g$ such that $g(\cdot, \cdot)=g(J\cdot, J\cdot)$. A manifold $M$ is called almost Hermitian manifold if it is endowed with an almost Hermitian structure $(J,g)$. An almost Hermitian manifold $(M, J, g)$ is called K\"{a}hler, if its Levi-Civita connection $\nabla$ satisfies $\nabla J=0$.

Recently, some mathematicians have studied some geometric concepts over Lie groups and Lie algebras such as complex, complex product  and contact structures \cite{ABD, AS, CF}. Inspired by these papers, the same authors in \cite{PN} introduced some geometric concepts  such as para-complex, para-Hermitian and para-K\"{a}hler structures on hom-Lie algebras. The aim of this paper is introducing other geometric concepts on these algebras.

The structure of this paper is organized as follows: In Section 2, we recall the definition of hom-algebra, hom-left-symmetric algebra, hom-Lie algebra, symplectic hom-Lie algebra and pseudo-Riemannian hom-algebra and hom-Levi-Civita product. Also we present an example of a symplectic hom-Lie algebra. In Section 3, we introduce complex and Hermitian structures on hom-Lie algebras and we determine all complex hom-Lie algebra of dimension 2. Also, we show that there not exists a proper complex (Hermitian) hom-Lie algebra of dimension 2. Moreover, we present an example of an almost Hermitian hom-Lie algebra. Then we study the complexification of a hom-Lie algebra. In Section 4, we define the phase space of a hom-Lie algebra and then using an involutive hom-left symmetric algebra we construct a phase space. Also, we present a complex structure on this phase space. In Section 5, we introduce  K\"{a}hler home Lie algebras and we present an example of them. Then we study some properties of them and we show that using these structures we can construct phase spaces.





\section{ hom-algebras and pseudo-Riemannian metric on hom-Lie algebra}
In this section, we present the definition of hom-algebra, hom-left-symmetric algebra, hom-Lie algebra and hom-Lie subalgebra. Then, we introduce symplectic hom-Lie algebra and pseudo-Riemannian hom-algebra.
\begin{definition}
A hom-algebra is a triple $(V,\cdot , \phi)$  consisting of a linear space $V$, a bilinear map (product)
$\cdot : V \times V \rightarrow V$ and an algebra morphism $\phi : V \rightarrow V$.
\end{definition}

Let $(V, \cdot, {\phi})$ be a hom-algebra. Then for any $u \in V$, we consider maps $L_u, R_u : V \rightarrow V$ as the left and the right
multiplication by u given by $L_u(v) = u\cdot v$ and $ R_u(v) = v\cdot u$, respectively. The ‎‎\textit{commutator}
 on $V$ is given by $[u,v]=u\cdot v-v\cdot u$.
 \begin{definition}
  Let $(V, \cdot, \phi)$ be a hom-algebra and $[\cdot,\cdot]$ be its commutator. We define the tensor curvature $\mathcal{K}$ of $V$ as follows
 \begin{equation}\label{SS8}
 \mathcal{K}(u,v):=L_{\phi(u)}\circ L_v-L_{ \phi(v)}\circ L_u-L_{[u,v]}\circ {\phi},
 \end{equation}
  for any $u, v\in V$.
 \end{definition}
 \begin{lemma}\label{IMPORT}
 We have 
\begin{equation}\label{SS6}
\circlearrowleft_{u,v,w}[{ \phi}(u),[v,w]]=\circlearrowleft_{u,v,w} \mathcal{K}(u,v)w,
\end{equation}
where $\circlearrowleft_{u,v,w}$ denotes cyclic sum on $u,v,w$. This equation is called hom-Bianchi identity.
 \end{lemma}
\begin{definition}
A hom-left-symmetric algebra is a hom-algebra $(V, \cdot,  \phi )$ such that
\begin{equation*}
ass_\phi(u,v,w)=ass_\phi(v,u,w),
\end{equation*}
where
\begin{align*}
ass_\phi(u,v,w)=(u\cdot v)\cdot\phi(w)-\phi(u)\cdot(v\cdot w),
\end{align*}
for any $u,v,w\in V$.
\end{definition}
\begin{definition}
A hom-Lie algebra is a triple $(\mathfrak{g}, [\cdot  , \cdot ],{ \phi‎‎_\mathfrak{g}})$ consisting of a linear space $\mathfrak{g}$, a
 bilinear map (bracket) $[\cdot , \cdot ]: \mathfrak{g}\times\mathfrak{g}\rightarrow\mathfrak{g}$ and an algebra morphism ${ \phi‎‎_\mathfrak{g}}:\mathfrak{g}\rightarrow \mathfrak{g}$ satisfying
\begin{align*}
[u,v]=-[v,u],
\end{align*}
\begin{equation}\label{AM502}
\circlearrowleft_{u,v,w}[\phi(u),[v,w]]=0,
\end{equation}
for any $u,v,w\in\mathfrak{g}$. The second equation is called hom-Jacobi identity. The hom-Lie algebra $(\mathfrak{g}, [\cdot  ,\cdot ], { \phi_\mathfrak{g}})$ is called regular (involutive), if ${ \phi_\mathfrak{g}}$ is non-degenerate (satisfies ${ \phi_\mathfrak{g}}^2=1$).
\end{definition}
\begin{definition}
Let $(\mathfrak{g}, [\cdot  , \cdot ],{ \phi‎‎})$ be a hom-Lie algebra. A subspace $\mathfrak{h}\subset \mathfrak{g}$ is called a hom-Lie subalgebra of $\mathfrak{g}$ if 
$
\phi(\mathfrak{h})\subset \mathfrak{h}$ and 
$[u,v]\in\mathfrak{h}$, for any $ {u,v\in \mathfrak{h}}.
$
\end{definition}
\begin{definition}
 A hom-algebra $(V, \cdot, \phi)$ is called
hom-Lie-admissible algebra
if its commutator bracket 
satisfies the hom-Jacobi identity.
\end{definition}
From Lemma \ref{IMPORT} and the above definition, results the following
\begin{corollary}
For a hom-Lie-admissible algebra we have $\circlearrowleft_{u,v,w} \mathcal{K}(u,v)w=0$.
\end{corollary}
\begin{proposition}\label{3.9}
If $(V, \cdot, \phi)$ is a hom-Lie-admissible algebra, then $(V, [\cdot ,\cdot ], \phi)$ is a hom-Lie algebra, where $[\cdot ,\cdot ]$ is the commutator bracket. 
\end{proposition}
Easily we can conclude the following
\begin{proposition}\label{3.10}
A hom-left-symmetric algebra is a hom-Lie-admissible algebra.
\end{proposition}
\begin{definition}
A symplectic hom-Lie algebra is a regular hom-Lie algebra $(\mathfrak{g}, [\cdot,\cdot],\phi_{g})$ endowed with a bilinear skew-symmetric
nondegenerate form $\omega$ which is $2$-hom-cocycle, i.e.
\begin{align*}
 \omega([u, v], { \phi_\mathfrak{g}}(w)) &+ \omega([w,u],{ \phi_\mathfrak{g}}( v)) + \omega([v,w],{ \phi_\mathfrak{g}}( u))=0,\\
&\omega(\phi(u),\phi(v))=\omega(u,v).
\end{align*}
In this case, $\omega$ is called a symplectic structure on $\mathfrak{g}$ and $(\mathfrak{g}, [\cdot,\cdot], \phi_\mathfrak{g}, \omega)$ is called symplectic hom-Lie algebra.
\end{definition}
\begin{theorem}\cite{PN}
Let $(\mathfrak{g}, [\cdot, \cdot], \phi_{g}, \omega)$ be an involutive  symplectic hom-Lie algebra. Then there exists a
 hom-left-symmetric algebra structure $\bold{a}$ on $\mathfrak{g}$ satisfying
\begin{equation}\label{AM10}
\omega(\bold{a}(u, v), { \phi_\mathfrak{g}}(w)) = - \omega({ \phi_\mathfrak{g}}(v),[u, w] ),
\end{equation}
such that
\[
\bold{a}(u, v) -\bold{a}(v, u) = [u, v],
\]
for any $u,v\in \mathfrak{g}$.
\end{theorem}
\begin{example}\label{IMEX}
We consider a $4$-dimensional linear space $\mathfrak{g}$ with an arbitrary basis $\{e_1,e_2,e_3,e_4\}$. We define the bracket
and linear map $\phi$ on $\mathfrak{g}$ as follows 
\begin{align*}
[e_1,e_2]=-ae_3,\ \ \  [e_1,e_3]=be_2,\ \ \ \ [e_2,e_4]=-a e_2,\ \ \ [e_3,e_4]=&ae_3,\ \ \
\end{align*}
and
\[
\phi(e_1)=-e_1,\ \ \ \phi(e_2)=e_2,\ \ \ \phi(e_3)=-e_3, \ \ \phi(e_4)=e_4.
\]
The above bracket is not a Lie bracket on $\mathfrak{g}$ if $a\neq 0$ and $b\neq0$, because
\[
[e_1,[e_3,e_4]]+[e_3,[e_4,e_1]+[e_4,[e_1,e_3]]=[e_1,ae_3]+[e_4,be_2]=2abe_2.
\]
It is easy to see that 
\begin{align*}
[\phi(e_1),\phi(e_2)]=&ae_3=\phi([e_1, e_2]),\ \ \ \ \ \   \ [\phi(e_1),\phi(e_3)]=be_2=\phi([e_1, e_3]),\\
 [\phi(e_2),\phi(e_4)]=&-ae_2=\phi([e_2, e_4]),\ \ \ \ [\phi(e_3),\phi(e_4)]=-a e_3=\phi([e_3, e_4]),
\end{align*}
i.e., $\phi$ is the algebra morphism. Also, we can deduce 
\begin{equation*}
[\phi(e_i),[e_j,e_k]]+[\phi(e_j),[e_k,e_i]+[\phi(e_k),[e_i,e_j]]=0,\ \ \ i,j,k=1, 2, 3, 4.
\end{equation*}
Thus $(\mathfrak{g}, [\cdot,\cdot], \phi)$ is a hom-Lie algebra. Now we consider the bilinear skew-symmetric
nondegenerate form $\Omega$ as follows:
\begin{equation}
\begin{bmatrix}
0&0&-A&0\\
0&0&0&\frac{a}{b}A\\
A&0&0&0\\
0&-\frac{a}{b}A&0&0
\end{bmatrix}
.
\end{equation}
Then we get
\begin{align*}
\Omega(\phi(e_1),\phi(e_3))=&-A=\Omega(e_1,e_3), \ \ \ \Omega(\phi(e_2),\phi(e_4))=\frac{a}{b}A=\Omega(e_2,e_4),\\
\Omega(\phi(e_1),\phi(e_2))=&0=\Omega(e_1,e_2),\ \ \ \ \  \Omega(\phi(e_1),\phi(e_4))=0=\Omega(e_1,e_4),\\
\Omega(\phi(e_2),\phi(e_3))=&0=\Omega(e_2,e_3),\ \ \ \ \  \Omega(\phi(e_3),\phi(e_4))=0=\Omega(e_3,e_4),
\end{align*}
and 
\begin{align*}
\Omega([e_i,e_j],\phi(e_k))+\Omega([e_j,e_k],\phi(e_i))+\Omega([e_k,e_i],\phi(e_j))=0,\ \ \ \ \  i,j,k=1, 2, 3, 4.
\end{align*}
The above relations show that $\Omega$ is $2$-hom-cocycle and so $(\mathfrak{g}, [\cdot,\cdot], \phi,\Omega)$ is a symplectic hom-Lie algebra.
\end{example}
\begin{definition}
Let $ (\mathfrak{g}, [\cdot, \cdot], { \phi_\mathfrak{g}})$ be a finite-dimensional hom-Lie algebra endowed
with a bilinear symmetric non-degenerate form $< , >$ such that for any $u,v\in \mathfrak{g}$ the following equation is satisfied
\begin{equation}\label{AM302}
	<\phi_\mathfrak{g}(u), \phi_\mathfrak{g}(v)>=<u, v>.
\end{equation}
Then, we say that $\mathfrak{g}$ admits a pseudo-Riemannian metric $<,>$ and $(\mathfrak{g}, [\cdot,\cdot], \phi_\mathfrak{g}, <,>)$  is called pseudo-Riemannian hom-Lie algebra.
\end{definition}
\begin{corollary}\label{AM350}
	Let $(\mathfrak{g}, [\cdot, \cdot], \phi_\mathfrak{g}, <,>)$ be a pseudo-Riemannian involutive hom-Lie algebra. Then for any $u,v\in\mathfrak{g}$ we have 
	\begin{equation}\label{EL1}
	<\phi_\mathfrak{g}(u), v>=<u, \phi_\mathfrak{g}(v)>.
	\end{equation}
\end{corollary}
If we consider a pseudo-Riemannian hom-Lie algebra $(\mathfrak{g}, [\cdot, \cdot], \phi_\mathfrak{g}, <,>)$  such that $\phi_\mathfrak{g}$ is isomorphism, then exists the unique product $\cdot$ on it which is given by Koszul’s formula
\begin{equation}\label{Koszul}
	2<u\cdot v,{ \phi_\mathfrak{g}}(w)>=<[u,v],{ \phi_\mathfrak{g}}(w)>+<[w,v],{ \phi_\mathfrak{g}}(u)>+<[w,u],{ \phi_\mathfrak{g}}(v)>,
\end{equation}
 which satisfies
\begin{equation}\label{SL}
[u,v]=u\cdot v-v\cdot u,
\end{equation}
\begin{equation}\label{L2}
<u\cdot v,{ \phi_\mathfrak{g}}(w)>=-<{ \phi_\mathfrak{g}}(v),u\cdot w>.
\end{equation}
This product is called the hom-Levi-Civita (see \cite{PN} for more details).
\section{Complex structure on hom-Lie algebras}
In this section, we introduce (almost) complex and Hermitian structures on hom-Lie algebras and study almost complex structures on hom-Lie algebras of dimension $2$. Also we introduce  the complexification of a hom-Lie algebra. 
\begin{definition}
	An almost complex structure on an involutive hom-Lie algebra $(\mathfrak{g}, [\cdot  , \cdot ], { \phi_\mathfrak{g}})$, is an isomorphism $ J : ‎‎\mathfrak{g} \rightarrow ‎‎\mathfrak{g}$ such that
	$
	J^2 = -Id_‎‎{\mathfrak{g}}$ and ${ \phi_\mathfrak{g}}\circ J=J\circ{ \phi_\mathfrak{g}}$. 
Also, $(\mathfrak{g}, [\cdot  , \cdot ], { \phi_\mathfrak{g}}, J)$ is called almost complex hom-Lie algebra.\\
If the $m$-dimensional hom-Lie algebra $\mathfrak{g}$ admits an almost complex structure $J$, then
\[
(det(J))^2=det(J^2)=det(-Id_{\mathfrak{g}})=(-1)^m,
\]
which implies that $m$ is even.
\end{definition}
	The above equations deduce $({ \phi_\mathfrak{g}}\circ J)^2=-Id_‎‎{\mathfrak{g}}$. The \textit{Nijenhuis torsion} of ${{ \phi_\mathfrak{g}}\circ J}$ is defined by
\begin{equation}\label{w1}
	N_{{ \phi_\mathfrak{g}}\circ J}(u,v)=[({{ \phi_\mathfrak{g}}\circ J})u, ({{ \phi_\mathfrak{g}}\circ J})v]-{{ \phi_\mathfrak{g}}\circ J}[({{ \phi_\mathfrak{g}}\circ J})u, v]-{{ \phi_\mathfrak{g}}\circ J}[u,({{ \phi_\mathfrak{g}}\circ J})v]-[u,v],
\end{equation}
for all $u,v\in \mathfrak{g}$. An almost complex structure is called complex if $N_{{ \phi_\mathfrak{g}}\circ J}=0$. In the following for simplicity, we set $N=N_{{ \phi_\mathfrak{g}}\circ J}$.
\begin{definition}\label{L100}
	An almost Hermitian structure of a hom-Lie algebra $(\mathfrak{g}, [\cdot,\cdot],\phi_{\mathfrak{g}})$ is a pair $(J, <,>)$ consisting of an almost complex structure
	and a pseudo-Riemannian metric $<,>$, such that for each $u,v \in \mathfrak{g}$
	\begin{equation}\label{EL}
	<({ \phi_\mathfrak{g}}\circ J)u, ({ \phi_\mathfrak{g}}\circ J)v>= <u , v>.
	\end{equation}
Also the pair $(J, <,>)$ is called Hermitian structure if $N=0$. In this case, $(\mathfrak{g}, [\cdot, \cdot], \phi_\mathfrak{g}, J, <,>)$ is called Hermitian hom-Lie algebra.
\end{definition}
\subsection{Almost complex structure on hom-Lie algebras of dimension 2}
In \cite{PN2}, the same authors determined all  non-abelian hom-Lie algebras of dimension 2 and non-abelian pseudo-Riemannain hom-Lie algebra of dimension 2. All non-abelian hom-Lie algebras of dimension 2  are as $(\mathfrak{g}, [\cdot, \cdot], \widehat{\phi})$, $(\mathfrak{g}, [\cdot, \cdot], \overline{\phi})$ and $(\mathfrak{g}, [\cdot,\cdot ], \widetilde{\phi})$, where
\begin{equation}\label{BB1}
(\widehat{\phi}(e_1)=e_1, \ \widehat{\phi}(e_2)=e_2),\ \ (\overline{\phi}(e_1)=e_1, \ \overline{\phi}(e_2)=-e_2),\ \ (\widetilde{\phi}(e_1)=e_1+Be_2,\  \widetilde{\phi}(e_2)=-e_2,\ \ B\neq0).
\end{equation}
Also all non-abelian pseudo-Riemannain hom-Lie algebra of dimension 2 are as $(\mathfrak{g}, [\cdot, \cdot], \widehat{\phi}, <,>)$, $(\mathfrak{g}, [\cdot, \cdot], \overline{\phi}, \prec,\succ)$ and $(\mathfrak{g}, [, ], \widetilde{\phi}, \ll,\gg)$, where $\widehat{\phi}, \overline{\phi}$ and $\widetilde{\phi}$ are given by (\ref{BB1})
and $<,>$ is an arbitrary bilinear symmetric non-degenerate form and  $\prec,\succ$, $\ll,\gg$ are bilinear symmetric non-degenerate forms with the following conditions
\begin{align*}
[<,>]&=
\begin{bmatrix}
<e_1,e_1>&<e_1,e_2>\\
<e_1,e_2>&<e_2, e_2>
\end{bmatrix},\ \ \ <e_1,e_1><e_2,e_2>-<e_1,e_2>^2\neq 0,\\
[\prec,\succ]&=
\begin{bmatrix}
\prec e_1,e_1\succ&0\\
0&\prec e_2, e_2\succ
\end{bmatrix},\ \ \ \prec e_1,e_1\succ\neq 0,\ \prec e_2,e_2\succ\neq0,\\
[\ll,\gg]&=
\begin{bmatrix}
\ll e_1,e_1\gg&-\frac{B}{2}\ll e_2, e_2\gg\\
-\frac{B}{2}\ll e_2, e_2\gg&\ll e_2, e_2\gg
\end{bmatrix},\ \ \ \ll e_1,e_1\gg\ll e_2,e_2\gg-\frac{B^2}{4}\ll e_2,e_2\gg^2\neq0.
\end{align*}
\begin{proposition}\label{EProp5}\cite{PN2}
The hom-Levi-Civita product on the pseudo-Riemannain hom-Lie algebra  $(\mathfrak{g}, [\cdot, \cdot], \widehat{\phi}, <,>)$ is
\begin{equation}\label{LEI4}
e_1\cdot e_1=0,\ \ \ e_1\cdot e_2=0,\ \ \ e_2\cdot e_1=-e_2,\ \ \ e_2\cdot e_2=\frac{<e_2, e_2>}{<e_1, e_1>}e_1,
\end{equation}
if $<e_1, e_2>=0$, and
\begin{align}
e_1\cdot e_1&=\frac{<e_1, e_2>^2}{\det [<,>]}e_1-\frac{<e_1, e_1><e_1, e_2>}{\det [<,>]}e_2,\label{LEI5}\\
e_1\cdot e_2&=\frac{<e_2, e_2>}{\det [<,>]}e_1-\frac{<e_1, e_2>}{\det [<,>]}e_2,\label{LEI6}\\
e_2\cdot e_1&=\frac{<e_2, e_2>}{\det [<,>]}e_1-\frac{\det [<,>]+<e_1, e_2>}{\det [<,>]}e_2,\label{LEI7}\\
e_2\cdot e_2&=\frac{<e_2, e_2>^2}{\det [<,>]<e_1, e_2>}e_1-\frac{<e_2, e_2>}{\det [<,>]}e_2,\label{LEI8}
\end{align}
if $<e_1, e_2>\neq0$.
\end{proposition}
\begin{proposition}\label{EProp2}
All non-abelian almost complex hom-Lie algebras of dimension 2 are as $(\mathfrak{g}, [\cdot, \cdot], \widehat{\phi}, \widehat{J})$, where $\widehat{\phi}$,  is given by (\ref{BB1})
and $\widehat{J}$ has the following matrix presentation:
\begin{equation*}
[\widehat{J}]=
\begin{bmatrix}
a&b\\
c&-a
\end{bmatrix},
\end{equation*}
where $a^2+bc=-1$.
\begin{equation*}
\end{equation*}
\end{proposition}
\begin{proof}
Let $(\mathfrak{g}, [\cdot, \cdot], \phi_\mathfrak{g}, J)$ be an almost product hom-Lie algebra. Conditions $J^2(e_1)=-e_1$ and $J^2(e_2)=-e_2$ give
\begin{align}\label{esi55}
(\rho^1_1)^2+\rho^2_1\rho^1_2=-1,\ \ \rho^2_1(\rho^1_1+\rho^2_2)=0,\ \  \rho^1_2(\rho^1_1+\rho^2_2)=0,\ \  (\rho^2_2)^2+\rho^2_1\rho^1_2=-1.
\end{align}
Now we consider two cases:

\textbf{Case 1.} $\rho_2^2=-\rho^1_1$.\\
In this case, the second and the third equations of (\ref{esi55}) are hold and the first and the fourth equations of (\ref{esi55}) reduce to $(\rho^1_1)^2+\rho^2_1\rho^1_2=-1$.

\textbf{Case 2.} $\rho_2^2\neq-\rho^1_1$.\\
In this case, the second and the third equations of (\ref{esi55}) give $\rho^2_1=\rho^1_2=0$, and so we deduce $(\rho^1_1)^2=(\rho^2_2)^2=-1$, and this is not possible.\\
Therefore $J$ has the following matrix presentation:
\begin{equation}\label{esi6}
\begin{bmatrix}
a&b\\
c&-a\\
\end{bmatrix},\ \ \ a^2+bc=-1.
\end{equation}
If we consider $\widehat{\phi}=Id$, then $(\mathfrak{g}, [\cdot,\cdot], \widehat{\phi}, \widehat{J})$ with $\widehat{J}$ given by (\ref{esi6}) is an almost complex (hom-)Lie algebra. Now, we consider $\overline{\phi}$. If the matrix presentation of $\overline{J}$ is as (\ref{esi6}), then condition $(\overline{J}\circ\overline{\phi})(e_1)=(\overline{\phi}\circ \overline{J})(e_1)$ implies $b=0$ and $(\overline{J}\circ\overline{\phi})(e_2)=(\overline{\phi}\circ \overline{J})(e_2)$ yields $c=0$. Consequently, we get the contradiction $a^2=-1$ and so $(\mathfrak{g}, [\cdot,\cdot], \overline{\phi}, \overline{J})$  can not be an almost complex hom-Lie algebra. Finally, we consider $\widetilde{\phi}_\mathfrak{g}$. If the matrix presentation of $\widetilde{J}$ is  (\ref{esi6}), then condition $(\widetilde{J}\circ\widetilde{\phi})(e_2)=(\widetilde{\phi}\circ \widetilde{J})(e_2)$ implies $c=0$ and consequently, $a^2=-1$. Therefore $\widetilde{J}$ can not be  an almost complex structure on $(\mathfrak{g}, [\cdot,\cdot], \widetilde{\phi})$.
\end{proof}
\begin{proposition}\label{EProp4}
All non-abelian almost Hermitian hom-Lie algebras of dimension 2 are as $(\mathfrak{g}, [\cdot, \cdot], \widehat{\phi}, \widehat{J}_i, <,>_i)$, $i=1,2$, where $\widehat{\phi}$ is given by (\ref{BB1})
and $\widehat{J}_i$ and  $<,>_i$ have the following matrix presentations:
\begin{equation}\label{esi10}
[\widehat{J}_1]=
\begin{bmatrix}
a&d\\
h&-a
\end{bmatrix},\ \
[\widehat{J}_2]=
\begin{bmatrix}
0&d\\
-\frac{1}{d}&0
\end{bmatrix},
\end{equation}
\begin{equation}
[<,>_1]=
\begin{bmatrix}
<e_1,e_1>_1&-\frac{a}{d}<e_1,e_1>_1\\
-\frac{a}{d}<e_1,e_1>_1&-\frac{h}{d}<e_1,e_1>_1
\end{bmatrix},\ \
[<,>_2]=
\begin{bmatrix}
<e_1,e_1>_2&0\\
0&\frac{1}{d^2}<e_1,e_1>_2
\end{bmatrix},
\end{equation}
where $a,d, h\neq 0$, $a^2+dh=-1$.
\end{proposition}
\begin{proof}
In Proposition \ref{EProp2}, it is shown that $(\mathfrak{g},
[\cdot, \cdot], \widehat{\phi})$ admit 
almost complex structure  $\widehat{J}$,
 i.e., $[\widehat{J}]=
\begin{bmatrix}
a&d\\
h&-a
\end{bmatrix}$. We consider the following cases \\
\textbf{Case 1.} $a\neq 0$.\\
In this case,  $<(\widehat{J}\circ\widehat{\phi})(e_1), (\widehat{J}\circ\widehat{\phi})(e_1)>=<e_1,e_1>$ implies
\begin{equation}\label{A1}
<e_1,e_2>=\frac{1-a^2}{2ad}<e_1,e_1>-\frac{d}{2a}<e_2,e_2>.
\end{equation}
 The condition $<(\widehat{J}\circ\widehat{\phi})(e_2), (\widehat{J}\circ\widehat{\phi})(e_2)>=<e_2,e_2>$ gives
  \begin{equation}\label{A2}
<e_1,e_2>=\frac{h}{2a}<e_1,e_1>+\frac{a^2-1}{2ah}<e_2,e_2>.
\end{equation}
From $<(\widehat{J}\circ\widehat{\phi})(e_1), (\widehat{J}\circ\widehat{\phi})(e_2)>=<e_1,e_2>$ we deduce
\begin{equation}\label{A3}
<e_1,e_2>=-\frac{a}{2d}<e_1,e_1>+\frac{a}{2h}<e_2,e_2>.
\end{equation}
 (\ref{A1}) and (\ref{A2}) imply
 \begin{equation}\label{A4}
<e_2,e_2>=-\frac{h}{d}<e_1,e_1>.
 \end{equation}
Setting (\ref{A4}) in (\ref{A3}), we obtain $<e_1,e_2>=-\frac{a}{d}<e_1,e_1>$. We denote these structures in the assertion with index $1$ in the below.\\
\textbf{Case 2.} $a= 0$.\\
 In the case, $hd=-1$ and consequently $[\widehat{J}]=
\begin{bmatrix}
0&d\\
-\frac{1}{d}&0
\end{bmatrix}$. Similar calculates give
$<e_2,e_2>=\frac{1}{d^2}<e_1,e_1>$ and $<e_1,e_2>=0$. We denote these structures in the assertion with index $2$ in the below.
\end{proof}
\begin{corollary}\label{FIN}
There not exists non-abelian almost Hermitian proper hom-Lie algebra of dimension 2.
\end{corollary}
According to the above result,
 we present an example of a $4$-dimensional space.
\begin{example}
We consider a $4$-dimensional hom-Lie algebra $(\mathfrak{g},[\cdot,\cdot],\phi)$ with an arbitrary basis $\{e_1,e_2,e_3,e_4\}$ such that bracket $[\cdot, \cdot]$ and linear map $\phi$ on $\mathfrak{g}$ are defined as follows 
\begin{align*}
[e_1,e_3]=ae_1+ae_2,\ \ \  [e_2,e_4]=ae_1+ae_2,\ \ \ \ [e_3,e_4]=-a e_3+ae_4,
\end{align*}
and
\[
\phi(e_1)=e_2,\ \ \ \phi(e_2)=e_1,\ \ \ \phi(e_3)=e_4, \ \ \phi(e_4)=e_3.
\]
The above bracket is not a Lie bracket on $\mathfrak{g}$ if $a\neq 0$. Now we consider the metric $<,>$ of $\mathfrak{g}$ as follows:
\begin{equation}
\begin{bmatrix}
1&0&0&0\\
0&1&0&0\\
0&0&1&0\\
0&0&0&1
\end{bmatrix}
.
\end{equation}
It is easy to see that  $<\phi(e_i),e_j>=<e_i,\phi(e_j)>$, for all $i,j=1,2,3,4$.
 If isomorphism $J$ is determined as
\[
J(e_1)=e_4,\ \ \ J(e_2)=e_3,\ \ \ J(e_3)=-e_2, \ \ J(e_4)=- e_1,
\]
then we have 
\[
J^2(e_1)=J^2(e_2)=J^2(e_3)=J^2(e_4)=-Id_\mathfrak{g}.
\]
The above equations yield 
\begin{align*}
(J\circ\phi)e_1=&e_3=(\phi\circ J)e_1,\ \ \ (J\circ\phi)e_2=e_4=(\phi\circ J)e_2,\\
(J\circ\phi)e_3=&-e_1=(\phi\circ J)e_3,\ \ \ (J\circ\phi)e_4=-e_2=(\phi\circ J)e_4.
\end{align*}
Thus $J$ is an almost complex structure on $(g, [\cdot, \cdot], \phi)$.  Also it results that 
 $<({ \phi‎‎_\mathfrak{g}}\circ J)(e_i), ({ \phi‎‎_\mathfrak{g}}\circ J)(e_j)>=<e_i, e_j>$, for all $i, j=1,2,3,4$, and so $(\mathfrak{g}, [\cdot, \cdot], \phi, <,>, J)$ is an almost Hermitian hom-Lie algebra. 
\end{example}
\subsection{ complexification of a hom-Lie algebra}
We denote the complexification of a hom-Lie algebra $\mathfrak{g}$ by $\mathfrak{g}^{\mathbb{C}}$ and we define it as follows
\[
\mathfrak{g}^{\mathbb{C}}=\{u+iv\  |\  u,v\in\mathfrak{g}\}.
\]
If $J$ is an almost complex structure on $\mathfrak{g}$, we may extend $J$ and $\phi_{\mathfrak{g}}$ to $\mathfrak{g}^{\mathbb{C}}$ as follows
\[
J(u+iv)=Ju+iJv, \ \ \ \ \phi_{\mathfrak{g}}(u+iv)=\phi_{\mathfrak{g}}(u)+i\phi_{\mathfrak{g}}(v).
\]
Obviously $J^2=-\phi_{\mathfrak{g}}^2=(\phi_{\mathfrak{g}}\circ J)^2=-Id_{\mathfrak{g}^{\mathbb{C}}}$ and ${\mathfrak{g}}\circ J=J\circ{\mathfrak{g}}$ on $\mathfrak{g}^{\mathbb{C}}$. Now, we consider
\begin{align*}
\mathfrak{g}^{1,0}=&\{w\in\mathfrak{g}^{\mathbb{C}}\ |\ (\phi_{\mathfrak{g}}\circ J)w=iw\},\\
\mathfrak{g}^{0,1}=&\{w\in\mathfrak{g}^{\mathbb{C}}\ |\ (\phi_{\mathfrak{g}}\circ J)w=-iw\}.
\end{align*}
\begin{proposition}\label{F15}
Let  $\mathfrak{g}$ be a hom-Lie algebra with an almost complex structure $J$. Then\\
i)\ \ $\mathfrak{g}^{1,0}=\{ u-i(\phi_{\mathfrak{g}}\circ J)u\ |\ u\in\mathfrak{g}\}$ and $ \mathfrak{g}^{0,1}=\{ u+i(\phi_{\mathfrak{g}}\circ J)u\ |\ u\in\mathfrak{g}\}$,\\
ii)\ $\mathfrak{g}^{\mathbb{C}}=\mathfrak{g}^{1,0}\oplus\mathfrak{g}^{0,1}$,\\
iii)\ $\overline{\mathfrak{g}^{1,0}}=\mathfrak{g}^{0,1}.$
\end{proposition}
\begin{proof}
Let $w=u+iv\in\mathfrak{g}^{\mathbb{C}}$. Then $(\phi_{\mathfrak{g}}\circ J)w=iw$ if and only if
\[
(\phi_{\mathfrak{g}}\circ J)u+v+i((\phi_{\mathfrak{g}}\circ J)v-u)=0,
\]
which gives $v=-(\phi_{\mathfrak{g}}\circ J)u$. Similarly, we see that $(\phi_{\mathfrak{g}}\circ J)w=-iw$ if and only if $v=(\phi_{\mathfrak{g}}\circ J)u$. Thus we have $(i).$
Let $w=u+iv\in\mathfrak{g}^{\mathbb{C}}$. Setting 
\begin{align*}
w_1=&\frac{1}{2}\big(u-i(\phi_{\mathfrak{g}}\circ J)u+i(v-i(\phi_{\mathfrak{g}}\circ J)v)\big),\\
w_2=&\frac{1}{2}\big(u+i(\phi_{\mathfrak{g}}\circ J)u+i(v+i(\phi_{\mathfrak{g}}\circ J)v)\big),
\end{align*}
we get $w_1\in\mathfrak{g}^{1,0}$ and $ w_2\in\mathfrak{g}^{0,1}$. Also it is easy to  
 see that
$
w_1+w_2=u+iv=w
$. Therefore $\mathfrak{g}^{\mathbb{C}}=\mathfrak{g}^{1,0}+ \mathfrak{g}^{0,1}$. If  $w\in \mathfrak{g}^{1,0}\cap\mathfrak{g}^{0,1}‎$, i.e., $w=u-i(\phi_{\mathfrak{g}}\circ J)u$ and $w=v+i(\phi_{\mathfrak{g}}\circ J)v$, then we obtain
\[
u-i(\phi_{\mathfrak{g}}\circ J)u=v+i(\phi_{\mathfrak{g}}\circ J)v,
\] 
and so 
\[
v-u+i\big((\phi_{\mathfrak{g}}\circ J)v+(\phi_{\mathfrak{g}}\circ J)u\big)=0,
\] 
which gives  $u=v$ and $u=-v$. Hence $\mathfrak{g}^{1,0}\cap\mathfrak{g}^{0,1}‎=0$ and consequently $\mathfrak{g}^{\mathbb{C}}=\mathfrak{g}^{1,0}\oplus\mathfrak{g}^{0,1}$. Thus we have $(ii)$. Let $w=u+iv\in\mathfrak{g}^{1,0}$. Then we have $v=-(\phi_{\mathfrak{g}}\circ J)u$. Therefore we get 
\[
(\phi_{\mathfrak{g}}\circ J)(\overline w)=(\phi_{\mathfrak{g}}\circ J)(u-iv)=-(v+iu)=-i(u-iv)=-i(\overline w),
\]
which implies $(\overline w)\in\mathfrak{g}^{0,1}$. So we have $(iii)$.
\end{proof}
\begin{proposition}
Let $J$ be an almost complex structure on a hom-Lie algebra $(\mathfrak{g},[\cdot,\cdot],\phi)$. Then the following statements are
equivalent:\\
i)\ \ $\mathfrak{g}^{1,0}$ is a  hom-Lie subalgebra, i.e., $[\mathfrak{g}^{1,0}, \mathfrak{g}^{1,0}] ‎‎\subseteq‎\mathfrak{g}^{1,0}$ and $\phi(\mathfrak{g}^{1,0})\subset \mathfrak{g}^{1,0}$,\\
ii)\ \ $\mathfrak{g}^{0,1}$ is a hom-Lie subalgebra, i.e., $[\mathfrak{g}^{0,1}, \mathfrak{g}^{0,1}] ‎‎\subseteq‎\mathfrak{g}^{0,1}$ and $\phi(\mathfrak{g}^{0,1})\subset \mathfrak{g}^{0,1}$,\\
iii)\ $J$ is a complex structure on $\mathfrak{g}$.\\
\end{proposition}
\begin{proof}
Assume $(i)$ holds. Let  $\widetilde{z},\widetilde{w}\in\mathfrak{g}^{0,1}$. Then according to $(iii)$ of Proposition \ref{F15}, there exist $z,w\in\mathfrak{g}^{1,0} $ such that 
 ${\widetilde{z}}=\overline{z}$ and $\widetilde{w}=\overline{{w}}$. On the other hand
\begin{align*}
[z,w]=&[u-i(\phi_{\mathfrak{g}}\circ J)u, v-i(\phi_{\mathfrak{g}}\circ J)v]=[u,v]-[(\phi_{\mathfrak{g}}\circ J)u,(\phi_{\mathfrak{g}}\circ J)v]-i([u,(\phi_{\mathfrak{g}}\circ J)v]+[(\phi_{\mathfrak{g}}\circ J)u,v]),
\end{align*}
for $u,v\in\mathfrak{g}$.
Therefore we have
\begin{align*}
\overline{[z,w]}=&[u,v]-[(\phi_{\mathfrak{g}}\circ J)u,(\phi_{\mathfrak{g}}\circ J)v]+i([u,(\phi_{\mathfrak{g}}\circ J)v]+[(\phi_{\mathfrak{g}}\circ J)u,v])\\
=&[u+i(\phi_{\mathfrak{g}}\circ J)u, v+i(\phi_{\mathfrak{g}}\circ J)v]=[\overline{z},\overline{w}].
\end{align*}
According to $(i)$, since $[z,w]\in\mathfrak{g}^{1,0}$, then from Proposition \ref{F15} and the above equation, we conclude $\overline{[z,w]}\in\mathfrak{g}^{0,1}$.
Hence 
$
[\widetilde{z},\widetilde{w}]=[\overline{z},\overline{w}]=\overline{[z,w]}\in\mathfrak{g}^{0,1}$. Also, since $\phi_{\mathfrak{g}}(z)\in \mathfrak{g}^{1,0}$, then
$
\phi_{\mathfrak{g}}(\widetilde z)=\phi_{\mathfrak{g}}(\overline z)=\overline{\phi_{\mathfrak{g}}(z)}\in\mathfrak{g}^{0,1}.
$
Moreover $[\phi_{\mathfrak{g}}(z),\phi_{\mathfrak{g}}(w)]=\phi_{\mathfrak{g}}[z,w]$ implies
\begin{align*}
[\phi_{\mathfrak{g}}(\widetilde{z}),\phi_{\mathfrak{g}}(\widetilde{w})]=&
[\phi_{\mathfrak{g}}(\overline{z}),\phi_{\mathfrak{g}}(\overline{w})]=[\overline{\phi_{\mathfrak{g}}(z)},\overline{\phi_{\mathfrak{g}}(w)}]=\overline{[\phi_{\mathfrak{g}}(z),\phi_{\mathfrak{g}}(w)]}=\overline{\phi_{\mathfrak{g}}[z,w]}=\phi_{\mathfrak{g}}\overline{[z,w]}=\phi_{\mathfrak{g}}[\overline{z},\overline{w}]\\
&=\phi_{\mathfrak{g}}[\widetilde{z},\widetilde{w}].
\end{align*}
Therefore $(i)$ implies $(ii)$. In the similar way, $(ii)$ implies $(i)$ and so $(i)$ is equivalent to $(ii)$.
Now, let $u,v\in\mathfrak{g}$ and $w:=[u+i(\phi_{\mathfrak{g}}\circ J)u, v+i(\phi_{\mathfrak{g}}\circ J)v]$. Then we get
\begin{align*}
w=[u,v]-[(\phi_{\mathfrak{g}}\circ J)u,(\phi_{\mathfrak{g}}\circ J)v]+i([u,(\phi_{\mathfrak{g}}\circ J)v]+[(\phi_{\mathfrak{g}}\circ J)u,v]).
\end{align*}
The above equation implies
\begin{align*}
(\phi_{\mathfrak{g}}\circ J)w=(\phi_{\mathfrak{g}}\circ J)[u,v]-(\phi_{\mathfrak{g}}\circ J)[(\phi_{\mathfrak{g}}\circ J)u,(\phi_{\mathfrak{g}}\circ J)v]+i(\phi_{\mathfrak{g}}\circ J)([u,(\phi_{\mathfrak{g}}\circ J)v]+[(\phi_{\mathfrak{g}}\circ J)u,v]),
\end{align*}
and 
\begin{align*}
iw=i[u,v]-i[(\phi_{\mathfrak{g}}\circ J)u,(\phi_{\mathfrak{g}}\circ J)v]-([u,(\phi_{\mathfrak{g}}\circ J)v]+[(\phi_{\mathfrak{g}}\circ J)u,v]).
\end{align*}
Summing two last equations, we obtain
\[
(\phi_{\mathfrak{g}}\circ J)w+iw=-(\phi_{\mathfrak{g}}\circ J)N(u,v)-i(\phi_{\mathfrak{g}}\circ J)N(u,v).
\]
Therefore $w\in\mathfrak{g}^{0,1}$, i.e., $(\phi_{\mathfrak{g}}\circ J)w=-iw$ if and only if $N=0$, which
means that $(ii)$ is is equivalent to $(iii)$.
\end{proof}

\section{Phase spaces and Representation }
For any $u\in V$, we consider $u^*\in V^*$ given by $u^*(v)=<u, v>$, for all $v\in V$. The map $u\rightarrow u^*$ define a canonical isomorphism between $V$ and $V^*$. In this paper, we denote $u^*(v)$ by $\prec u^*, v\succ$.
\begin{definition}
Let $(V, \cdot , \phi)$ be a hom-algebra and $V^*$ be its dual space. If there exists a hom-algebra structure on the direct sum of the underling vector space $(V, \phi)$ and $(V^*, \phi^*)$ ($\phi^*$ is the dual map of $\phi$, i.e., $\phi^*=\phi^t$) such that $(V,\phi)$ and $(V^*,\phi^*)$ are sub-hom-algebras and the natural skew-symmetry bilinear form $\omega$
on $V\oplus V^*$ given by 
\begin{align}\label{S1}
\omega(u+a^* ,v+ b^* )=\prec b^* ,u‎\succ-\prec a^* ,v‎\succ,
\end{align}
is a symplectic form for any $u,v\in V$ and $a^* , b^* \in V^*$, then $(V\oplus V^*, \omega)$ is called phase space of $V$.
\end{definition}
\begin{definition}
Let $(\mathfrak{g}, [.,.],\phi_{\mathfrak{g}})$ be a hom-Lie algebra and $\mathfrak{g}^*$ be
its dual space. A phase space of $\mathfrak{g}$ is defined as a hom-Lie algebra $T^*\mathfrak{g}=(\mathfrak{g}\oplus‎\mathfrak{g}^*, [. ,.]_{\mathfrak{g}\oplus‎\mathfrak{g}^*}, \phi_{\mathfrak{g}}\oplus\phi_{\mathfrak{g^*}})$ such that it endowed with the symplectic form $\omega$ given by (\ref{S1}), where $\phi_{\mathfrak{g}^*}=(\phi_\mathfrak{g})^*$.
\end{definition}
Let $(g, [\cdot  , \cdot ], {\phi_\mathfrak{g}})$ be a hom-Lie algebra. A \textit{representation} of $\mathfrak{g}$ is a triple $(V, A, \rho)$ where $V$ is a vector space, $A\in gl(V)$ and 
$\rho: \mathfrak{g}\rightarrow gl(V)$ is a linear map satisfying
\begin{equation}\label{SS4}
\left\{
\begin{array}{cc}
&\hspace{-4cm}\rho(\phi_{\mathfrak{g}}(u))\circ A=A\circ \rho(u),\\
&\hspace{-.2cm}\rho([u,v]_{\mathfrak{g}})\circ A=\rho(\phi_{\mathfrak{g}}(u))\circ \rho(v)-\rho(\phi_{\mathfrak{g}}(v))\circ \rho(u),
\end{array}
\right.
\end{equation}
for any $u,v \in \mathfrak{g}$. If we consider $V^*$ as the dual vector
space of $V$ then we can define a linear map $\widetilde{\rho}:\mathfrak{g}\rightarrow gl(V^*)$ by 
\[
\prec\widetilde{\rho}(u)(a^* ),v\succ=-\prec({\rho}(u))^t(a^* ),v\succ=-\prec a^* ,\rho(u)(v)\succ,
\]
for any $u\in \mathfrak{g}, v\in V,  a^* \in {V^*}$, where $(\rho(u))^t$ is the dual endomorphism of $\rho(u)\in\End(V)$ and $\prec\widetilde{\rho}(u)(a^* ),v\succ$ is $\widetilde{\rho}(u)(a^* )(v)$. The representation $(V,A, \rho)$ is called \textit{admissible} if $(V ^*,A^*=A^t,\widetilde{\rho})$ is also a representation of $\mathfrak{g}$. In \cite{SB}, it is shown that the representation $(V, A, {\rho})$ is admissible if and only if 
the following conditions  are satisfied  
\begin{equation}\label{AM200}
\left\{
\begin{array}{cc}
&\hspace{-3.7cm}A\circ \rho(\phi_{\mathfrak{g}}(u))=\rho(u)\circ A,\\
&A\circ\rho([u,v]_\mathfrak{g})=\rho(u)\circ \rho(\phi_{\mathfrak{g}}(v))-\rho(v)\circ\rho(\phi_{\mathfrak{g}}(u)).
\end{array}
\right.
\end{equation}
Let $(g, [\cdot  , \cdot ], {\phi_\mathfrak{g}})$ be a hom-Lie algebra and $ad : \mathfrak{g}\rightarrow End(\mathfrak{g})$ be an operator defined for any
$u,v\in \mathfrak{g}$ by $ad(u)(v) = [u, v]$.
Moreover, it is easy to see that  
\[
ad{[u,v]}\circ{ \phi_\mathfrak{g}}=ad({ \phi_\mathfrak{g}}(u))\circ ad(v)-ad({ \phi_\mathfrak{g}}(v))ad(u).
\]
Thus $(\mathfrak{g}, ad, { \phi_\mathfrak{g}})$ is a representation of $\mathfrak{g}$, which is called an \textit{adjoint representation of $\mathfrak{g}$} (see \cite{SB}, for more details).
\begin{lemma}\cite{PN}
Let $(g, [\cdot  , \cdot ], {\phi_\mathfrak{g}})$ be an involutive hom-Lie algebra. Then the adjoint representation  $(\mathfrak{g}, { \phi_\mathfrak{g}},ad)$ is
admissible.
\end{lemma} 
\begin{proposition}
Let $(V, \cdot,\phi)$ be a  hom-left symmetric algebra and for
  any $u\in V$,  $L_u$ be the left
multiplication operator on $u\in V$ (i.e., $L_uv=u\cdot v$ for any $v\in V)$. Then 
$(V,\phi,L)$ is a representation of the hom-Lie algebra $V$,  where
$L: V\rightarrow gl(V)$ with $u\rightarrow L_u$.
\end{proposition}
\begin{proof}
We must to show that 
\begin{equation}\label{5.2}
\left\{
\begin{array}{cc}
\ \ \ \ \ \ \ \  \ \ \ \ \ \ \ \ &\hspace{-5cm}L_{ \phi(u)}\circ{ \phi‎‎}={ \phi‎‎}\circ L_u,\\
\ \ \ \ \ \ \ \  \ \ &\hspace{-2.3cm}L_{[u,v]}\circ { \phi}=L_{ \phi‎‎(u)}\circ L_v-L_{ \phi‎‎‎(v)}\circ L_u.
\end{array}
\right.
\end{equation}
If we consider the hom-left-symmetric algebra $V$, then for any $u,v,w\in V$ we have
\begin{align*}
(u\cdot v)\cdot{ \phi‎‎}(w)-{ \phi‎‎}(u)\cdot (v\cdot w)=(v\cdot u)\cdot{ \phi}(w)-{ \phi}(v)\cdot(u\cdot w).
\end{align*}
Setting $u\cdot v-v\cdot u=[u,v]$ and $L_uv=u\cdot v$, in the above equation we get
 \[
L_{[u,v]}\circ { \phi}=L_{{ \phi}(u)}\circ L_v-L_{ { \phi}‎(v)}\circ L_u.
\]
Moreover, we have ${ \phi}(u\cdot v)={ \phi}(u)\cdot { \phi}(v)$, that is
$L_{{ \phi}(u)}{ \phi}(v)={ \phi}(L_uv)$. Thus we have (\ref{5.2}).
\end{proof}
\begin{corollary}
Let $(V, \cdot, \phi)$ be an  involutive  hom-left symmetric algebra. Then representation $(V,\phi,L)$ is admissible, i.e.,  $(V^*,\phi^*,\widetilde{L})$ is a  representation of $V$. In particular, $\phi^*$ is  involutive.
\end{corollary}
\begin{proof}
According to (\ref{AM200}), it is sufficient to prove
\begin{equation*}
\left\{
\begin{array}{cc}
 &\hspace{-3cm}{ \phi}\circ L_{{ \phi}(u)}=L_u\circ { \phi},\\
&{ \phi}\circ L_{[u,v]}=L_u\circ L_{ \phi(v)}-L_v\circ L_{ \phi(u)}.
\end{array}
\right.
\end{equation*}
If $u,v\in V$, then using ${ \phi}(u\cdot v)={ \phi}(u)\cdot { \phi}(v)$, we obtain ${ \phi}({ \phi}(u)\cdot v)=u\cdot { \phi}(v)$. Considering ${ \phi}(u)\cdot v=L_{ \phi(u)}v$, we deduce ${ \phi}( L_{{ \phi}(u)}v)=L_u { \phi}(v)$.
On the other hand, for any $w\in V$ we have 
\[
(u\cdot v)\cdot{ \phi}(w)-(v\cdot u)\cdot{ \phi}(w)={ \phi}(u)\cdot(v\cdot w)-{ \phi}(v)\cdot(u\cdot w).
\]
Contracting the above equation with ${ \phi}$ implies that
\[
{ \phi}((u\cdot v)\cdot{ \phi}(w)-(v\cdot u)\cdot{ \phi}(w))=u\cdot{ \phi}(v\cdot w)-v\cdot{ \phi}(u\cdot w).
\]
Setting $L_uv=u\cdot v$ in the above equation yields
\[
{ \phi}(L_{[u,v]}{ \phi}(w))=L_uL_{ \phi(v)}{ \phi}(w)-L_vL_{ \phi(u)}{ \phi}(w).
\]
Therefore we have the assertion.
\end{proof}
\begin{proposition}
Let $(V,\cdot, \phi)$ be an involutive hom-left-symmetric algebra and $(V,\phi,L)$ be a representation of the hom-Lie algebra $V$. Then $(V\oplus V^*,\cdot,\Phi)$ is a hom-left symmetric algebra where $\cdot$ and $\Phi$ are given by 
\begin{equation}\label{F2}
\left\{
\begin{array}{cc}
 &\hspace{-.4cm}
(u,a^*)\cdot(v,b^*):=(u\cdot v,\widetilde{L}_{\phi(u)}b^*),\\
&\Phi(u,a^*):=(\phi(u),\phi^*(a^*)),
\end{array}
\right.
\end{equation}
for any $u,v\in V, a^*,b^*\in V^*$.
\end{proposition}
\begin{proof}
We have
\begin{align*}
\Phi((u,a^* )\cdot(v,b^* ))=\Phi(u\cdot v,\widetilde{L}_{\phi(u)}b^* )=(\phi(u\cdot v),\phi^*(\widetilde{L}_{\phi(u)}b^* ))=(\phi(u)\cdot \phi(v),\phi^*(\widetilde{L}_{\phi(u)}b^* )),
\end{align*}
and so
\begin{align*}
\Phi(u,a^* )\cdot\Phi(v,b^* )=(\phi(u),\phi^*(a^* ))\cdot(\phi(v),\phi^*(b^* ))=(\phi(u)\cdot\phi(v),\widetilde{L}_u\phi^*(b^* )).
\end{align*}
But $L$ is an admissible representation on $V$, i.e., $\phi^*(\widetilde{L}_{{\phi(u)}}b^* )=\widetilde{L}_u\phi^*(b^* )$. Thus 
\[
\Phi((u,a^* )\cdot(v,b^* ))=\Phi(u,a^* )\cdot\Phi(v,b^* ),
\]
which implies $\Phi$ is a home algebra structure on $V\oplus V^*$. Now, we show that the product given by (\ref{F2}) is hom-left symmetric. A direct computation yields
\begin{align}\label{NK2}
&((u,a^* )\cdot(v, b^* ))\cdot\Phi(w, c^*)-\Phi(u,a^* )\cdot((v, b^* )\cdot(w, c^*))-
((v, b^* )\cdot(u,a^* ))\cdot\Phi(w, c^*)\\
&+\Phi(v, b^* )\cdot((u,a^* )\cdot(w, c^*))\nonumber\\
&=(L_{[u,v]} { \phi}(w)-L_{ \phi‎‎(u)} L_vw+L_{ \phi‎‎‎(v)} L_uw,\widetilde{L}_{\phi([u,v])}{ \phi}^*( c^*)-\widetilde{L}_{ u} \widetilde{L}_{\phi(v)} c^*+\widetilde{L}_{ v} \widetilde{L}_{\phi(u)} c^*).\nonumber
\end{align}
On the other hand, we have
\begin{align}\label{NK1}
L_{[u,v]} { \phi}(w)=L_{ \phi‎‎(u)} L_vw-L_{ \phi‎‎‎(v)} L_uw,\nonumber\\
\widetilde{L}_{[u,v]}{ \phi}^*( c^*)=\widetilde{L}_{ \phi‎‎(u)} \widetilde{L}_v c^*-\widetilde{L}_{ \phi‎‎‎(v)} \widetilde{L}_u c^*.
\end{align}
 Putting $u:=\phi(u), v:=\phi(v)$ in (\ref{NK1}) we deduce
\begin{align*}
\widetilde{L}_{\phi([u,v])}{ \phi}^*( c^*)=\widetilde{L}_{ u} \widetilde{L}_{\phi(v)} c^*-\widetilde{L}_{ v} \widetilde{L}_{\phi(u)} c^*.
\end{align*}
Setting the above equations in (\ref{NK2})  we conclude the assertion. 
\end{proof}
\begin{proposition}
Hom-left symmetric algebra $(V\oplus V^*,\cdot,\Phi)$ is an involutive hom-Lie algebra with commutator bracket, where $\cdot$ and $\Phi$ are given by (\ref{F2}).
\end{proposition}
\begin{proof}
First, we show that $(V\oplus V^*,\cdot,\Phi)$ is a hom-Lie-admissible algebra. If we set
\begin{align*}
[(u,a^* ),(v,b^* )]=(u,a^* )\cdot(v,b^* )-(v,b^* )\cdot(u,a^* ),\ \ \forall u,v\in V, \ \ \forall a^* ,b^* \in V^*,
\end{align*}
then using (\ref{F2}) we have
\begin{align}\label{F4}
[(u,a^* ),(v,b^* )]&=(u\cdot v-v\cdot u,\widetilde{L}_{\phi(u)}b^* -\widetilde{L}_{\phi(v)} a^* ).
\end{align}
Obviously, we have 
\[
[(u,a^* ),(v,b^* )]=-[(v, b^* ),(u,a^* )].
\]
Now, we check the hom-Jacobi identity. We have
\begin{align*}
&\circlearrowleft‎‎‎_{(u,a^*  ),(v, b^*  ),(w, c^* )}[\Phi(u,a^*  ),[(v, b^*  ),(w, c^* )]]=\circlearrowleft‎‎‎_{(u,a^*  ),(v, b^*  ),(w, c^* )}[(\phi(u),\phi^*(a^*  )),(v\cdot w-w\cdot v,\widetilde{L}_{\phi(v)} c^* \\
&-\widetilde{L}_{\phi(w)}  b^*  )]
=\circlearrowleft‎‎‎_{(u,a^*),(v, b^*  ),(w, c^* )}(\phi(u)\cdot(v\cdot w-w\cdot v)-(v\cdot w-w\cdot v)\cdot\phi(u),\widetilde{L}_{u}( \widetilde{L}_{\phi(v)} c^* -\widetilde{L}_{\phi(w)}  b^*  )\\
&-\widetilde{L}_{\phi([v,w])}\phi^*(a^*  )).
\end{align*}
Since $V$ is a hom-left-symmetric algebra, then
\[
(u\cdot v)\cdot\phi(w)-\phi(u)\cdot(v\cdot w)-(v\cdot u)\cdot\phi(w)+\phi(v)\cdot(u\cdot w)=0.
\]
On the other hand, $(V^*,\phi^*,\widetilde{L})$ is a representation of $V$. Therefore
 \[
\widetilde{L}_{\phi([v,w])}\circ { \phi}^*=\widetilde{L}_v\circ \widetilde{L}_{\phi(w)}-\widetilde{L}_w\circ \widetilde{L}_{\phi(v)}.
\]
The above equations imply 
\begin{align*}
&\circlearrowleft‎‎‎_{(u,a^*  ),(v, b^*  ),(w, c^* )}[\Phi(u,a^*  ),[(v, b^*  ),(w, c^* )]]=(\phi(u)\cdot(v\cdot w-w\cdot v)-(v\cdot w-w\cdot v)\cdot\phi(u)\\
&+
\phi(v)\cdot(w\cdot u-u\cdot w)
-(w\cdot u-u\cdot w)\cdot\phi(v)
+\phi(w)\cdot(u\cdot v-v\cdot u)-(u\cdot v-v\cdot u)\cdot\phi(w)\\
&,\widetilde{L}_u( \widetilde{L}_{{ \phi}(v)}  c^* -\widetilde{L}_{{ \phi}(w)}  b^*  )-(\widetilde{L}_v\circ \widetilde{L}_{{ \phi}(w)}-\widetilde{L}_w\circ \widetilde{L}_{{ \phi}(v)})(a^*  )+
\widetilde{L}_v( \widetilde{L}_{{ \phi}(w)} a^*  -\widetilde{L}_{{ \phi}(u)}  c^* )\\
&-(\widetilde{L}_w\circ \widetilde{L}_{{ \phi}(u)}-\widetilde{L}_u\circ \widetilde{L}_{{ \phi}(w)})( b^*  )
+\widetilde{L}_w(\widetilde{L}_{{ \phi}(u)} b^*  -\widetilde{L}_{{ \phi}(v)}a^*  )-(\widetilde{L}_u\circ \widetilde{L}_{{ \phi}(v)}-\widetilde{L}_v\circ \widetilde{L}_{{ \phi}(u)})( c^* )=0,
\end{align*}
which conclude the hom-Jacobi identity. Using (\ref{F2}) we obtain
\begin{align*}
\Phi[(u,a^* ),(v,b^* )]=&\Phi(u\cdot v-v\cdot u,\widetilde{L}_{\phi(u)}b^* -\widetilde{L}_{\phi(v)}a^* )=(\phi(u\cdot v-v\cdot u),\phi^*(\widetilde{L}_{\phi(u)}b^* -\widetilde{L}_{\phi(v)}a^*  ))\\
=&(\phi(u)\cdot\phi( v)-\phi(v)\cdot\phi( u),\widetilde{L}_u\phi^*(b^* )-\widetilde{L}_v\phi^*(a^* ))\\
=&[(\phi(u),\phi^*( a^* )),({\phi(v)},\phi^*(b^* ))]=[\Phi(u,a^* ),\Phi(v,b^* )].
\end{align*}
Thus $(V\oplus V^*,\cdot,\Phi)$ is a hom-Lie algebra with commutator bracket. Finally, we get 
\begin{align*}
\Phi^2(u,a^* )=\Phi(\phi(u),\phi^*(a^* ))=(\phi^2(u),(\phi^*)^2(a^* ))=(u,a^*),
\end{align*}
i.e., $\Phi$ is involutive.
\end{proof}
\begin{proposition}
Let $(V,\cdot, \phi)$ be an involutive hom-left-symmetric algebra and $(V,\phi,L)$ be a representation of the hom-Lie algebra $V$. Then
\begin{align}\label{AM600}
‎‎\prec \widetilde{L}_{\phi(u)}a^*,\phi(v)\succ =-‎‎\prec u\cdot v,\phi^*( a^*)‎\succ.
\end{align}
\end{proposition}\label{F11}
\begin{proof}
We have
\begin{align*}
\prec\phi^*( a^*), u\cdot v‎\succ‎‎&=‎\prec  a^*,\phi(u\cdot v)‎\succ=‎‎\prec  a^*,\phi(u)\cdot\phi( v)‎\succ
‎=\prec a^*,L_{\phi(u)}\phi( v)‎\succ‎\\&=‎‎\prec L_{\phi(u)}^t a^*,\phi(v) ‎\succ‎=‎-‎‎\prec \widetilde{L}_{\phi(u)}a^*,\phi(v) ‎\succ,
\end{align*}
for any $u,v\in V$ and $ a^*\in V^*$. 
\end{proof}
\begin{theorem}
Let $(V\oplus V^*,\Phi, [\cdot,\cdot])$ be a hom-Lie algebra where $\Phi$ and $[\cdot,\cdot]$ are given by (\ref{F2}) and (\ref{F4}). If we consider a non-degenerate
bilinear form $\Omega_{_{V\oplus V^*}}$ on $V\oplus V^*$ given by
\begin{equation}\label{AM505}
\Omega_{_{V\oplus V^*}}(‎ (u,a^* ), (v, b^*))=‎‎\prec‎  b^*,u\succ-‎‎\prec a^*,v‎\succ,\ \ \ \ \ \forall u,v\in V, \ \ \ \forall  a^*, b^*\in V^*,
\end{equation}
then $V\oplus V^*$ is a phase space.
\end{theorem}
\begin{proof}
Using (\ref{F2}) and (\ref{AM505}), we get
\begin{align*}
\Omega_{_{V\oplus V^*}}(‎ \Phi(u, a^*),\Phi( v, b^*))=&\Omega_{_{V\oplus V^*}}(‎ (\phi(u),\phi^*( a^*)),(\phi( v),\phi^*( b^*)))=‎‎\prec‎ \phi^*( b^*),‎ \phi(u)\succ-‎‎\prec\phi^*( a^*),\phi( v)‎\succ\\
&=‎‎\prec‎  b^*,u\succ-‎‎\prec a^*,v‎\succ=\Omega_{_{V\oplus V^*}}(‎ u+ a^*, v+ b^*).
\end{align*}
Also (\ref{AM600}) implies
\begin{align*}
&\circlearrowleft‎‎‎_{(u, a^*),(v, b^*),(w, c^*)}\Omega_{_{V\oplus V^*}}\big([‎ ‎ (u, a^*), (v, b^*)],\Phi((w, c^*))\big)\\
&=\circlearrowleft‎‎‎_{(u, a^*),(v, b^*),(w, c^*)}\Omega_{_{V\oplus V^*}}\big(u\cdot v-v\cdot u,\widetilde{L}_{\phi(u)}b^*-\widetilde{L}_{\phi(v)}a^*),(\phi(w),\phi^*(c^*))\big)=‎‎\prec‎  \phi^*(c^*),u\cdot v-v\cdot u\succ\\
&-‎‎\prec‎ \widetilde{L}_{\phi(u)}b^*-\widetilde{L}_{\phi(v)}a^*),\phi(w)\succ
+‎‎\prec‎  \phi^*(a^*),v\cdot w-w\cdot v\succ
-‎‎\prec‎ \widetilde{L}_{\phi(v)}c^*-\widetilde{L}_{\phi(w)}b^*),\phi(u)\succ\\
&+‎‎\prec‎  \phi^*(b^*),w\cdot u-u\cdot w\succ
-‎‎\prec‎ \widetilde{L}_{\phi(w)}a^*-\widetilde{L}_{\phi(u)}c^*),\phi(v)\succ=0.
\end{align*}
Thus $\Omega_{_{V\oplus V^*}}$ is a $2$-hom-cocycle form and consequently $V\oplus V^*$ is a phase space.
\end{proof}
\begin{lemma}\label{F12}
Let $(V, \cdot, \phi)$ be an  involutive  hom-left symmetric algebra. Then for representation $(V,\phi,L)$ we have
\[
\phi^*(u^*)=(\phi(u))^*.
\]
\end{lemma}
\begin{proof}
For any $v\in V$ we have
\begin{align*}
\phi^*(u^*)(v)=u^*(\phi(v))=<u, \phi(v)>,
\end{align*}
and 
\begin{align*}
(\phi(u))^*(v)=<\phi(u), v>.
\end{align*}
Since $\phi$ is involutive, then (\ref{EL1}) implies that the above equations are equal.
\end{proof}
\begin{proposition}\label{Prop4.11}
Let $(V\oplus V^*,\Phi,[\cdot,\cdot])$  be a  hom-Lie algebra where $\Phi$ and $[\cdot,\cdot]$ are given by (\ref{F2}) and (\ref{F4}). The linear map
$\mathcal{J}:V\oplus V^*\rightarrow V\oplus V^*$ defined by 
\begin{align}\label{F14}
\mathcal{J}(u,a^*)=(-\phi(a),\phi^*(u^*)),\ \ \ \forall u\in V, a^*\in V^*,
\end{align}
 is a complex structure on the hom-Lie algebra $V\oplus V^*$.
\end{proposition}
\begin{proof}
Applying Lemma \ref{F12}, we obtain 
\begin{align*}
\mathcal{J}^2(u,a^*)=&\mathcal{J}(-\phi(a),\phi^*(u^*))=\mathcal{J}(-\phi(a),(\phi(u))^*)=(-\phi(\phi(u)),-\phi^*(\phi(a))^*)=-(u,a^*).
\end{align*}
Also $\Phi\circ \mathcal{J}=\mathcal{J}\circ\Phi$, because
\begin{align*}
(\Phi\circ \mathcal{J})(u,a^*)=&\Phi(-\phi(a),\phi^*(u^*))=(-\phi^2(a),(\phi^*)^2(u^*))=(-a,u^*),\\
(\mathcal{J}\circ\Phi)(u,a^*)=&\mathcal{J}(\phi(u),\phi^*(a^*))=(-\phi(\phi(a)),\phi^*(\phi(u))^*)=(-a,u^*).
\end{align*}
Using (\ref{F2}) we have
\begin{align}\label{LLL}
&\mathcal{N}((u,a^*),(v,b^*))=[(\Phi\circ \mathcal{J})(u,a^*),(\Phi\circ \mathcal{J})(v,b^*)]-(\Phi\circ \mathcal{J})[(\Phi\circ \mathcal{J})(u,a^*),(v,b^*)]\\
&-(\Phi\circ \mathcal{J})[(u,a^*),(\Phi\circ \mathcal{J})(v,b^*)]-[(u,a^*),(v,b^*)]=[(-a,u^*),(-b,v^*)]-(\Phi\circ \mathcal{J})[(-a,u^*),(v,b^*)]\nonumber\\
&-(\Phi\circ \mathcal{J})[(u,a^*),(-b,v^*)]-[(u,a^*),(v,b^*)]=(a\cdot b-b\cdot a,-\widetilde{L}_{\phi(a)}v^*+\widetilde{L}_{\phi(b)}u^*)-(\Phi\circ \mathcal{J})(-a\cdot v\nonumber\\
&+v\cdot a,-\widetilde{L}_{\phi(a)}b^*-\widetilde{L}_{\phi(v)}u^*)-(\Phi\circ \mathcal{J})(-u\cdot b+b\cdot u,\widetilde{L}_{\phi(u)}v^*+\widetilde{L}_{\phi(b)}a^*)-(u\cdot v-v\cdot u,\widetilde{L}_{\phi(u)}b^*-\widetilde{L}_{\phi(v)}a^*)\nonumber.
\end{align}
Now let $u\in \mathfrak{g}$. Then we have 
\begin{align*}
\phi^*(\widetilde{L}_{\phi(a)}v^*)(u)=-\phi^*({L}^*_{\phi(a)}v^*)(u)=-{L}^*_{\phi(a)}v^*(\phi(u))=-v^*(L_{\phi(a)}\phi(u))=-<v,\phi(a)\cdot \phi(u)>.
\end{align*}
Applying (\ref{L2}) in the above equation we conclude 
\begin{align*}
\phi^*(\widetilde{L}_{\phi(a)}v^*)(u)=<\phi(a)\cdot\phi(v), u>=<L_{\phi(a)}\phi(v),u>=(L_{\phi(a)}\phi(v))^*(u),
\end{align*}
i.e., 
\begin{align*}
\phi^*(\widetilde{L}_{\phi(a)}v^*)=(L_{\phi(a)}\phi(v))^*.
\end{align*}
The last equation and $(\Phi\circ \mathcal{J})(u,a^*)=(-a,u^*)$ imply
\begin{align*}
&(\Phi\circ \mathcal{J})(a\cdot v,\widetilde{L}_{\phi(a)}b^*)=(\mathcal{J}\circ\Phi )(a\cdot v,\widetilde{L}_{\phi(a)}b^*)=\mathcal{J}(\phi(a)\cdot\phi( v),\phi^*(\widetilde{L}_{\phi(a)}b^*))=\mathcal{J}(L_{\phi(a)}\phi( v),\phi^*(\widetilde{L}_{\phi(a)}b^*))\\
=&\mathcal{J}((\phi^*(\widetilde{L}_{\phi(a)}v^*))^*,(L_{\phi(a)}\phi(b))^*)=(-\phi(L_{\phi(a)}\phi(b)),\phi^*(\phi^*(\widetilde{L}_{\phi(a)}v^*)))=(-a\cdot b,\widetilde{L}_{\phi(a)}v^*).
\end{align*}
Finally the above equation and $(\ref{LLL})$ yield
\begin{align*}
&\mathcal{N}((u,a^*),(v,b^*))=
(a\cdot b-b\cdot a-u\cdot v+v\cdot u-a\cdot b-v\cdot u+u\cdot v+b\cdot a,-\widetilde{L}_{\phi(a)}v^*+\widetilde{L}_{\phi(b)}u^*-\widetilde{L}_{\phi(u)}b^*\\
&+\widetilde{L}_{\phi(v)}a^*+\widetilde{L}_{\phi(a)}v^*-\widetilde{L}_{\phi(v)}a^*-\widetilde{L}_{\phi(b)}u^*+\widetilde{L}_{\phi(u)}b^*)=0.
\end{align*}
\end{proof}
\section{ K\"{a}hler hom-Lie algebras}
In this section, we introduce K\"{a}hler structures on hom-Lie algebras. Also, we present an example of these structures. 
\begin{definition}
A K\"{a}hler hom-Lie algebra is a pseudo-Riemannian hom-Lie algebra $(\mathfrak{g}, [\cdot, \cdot], \phi_{\mathfrak{g}},<,>)$ endowed with an almost complex structure $J$, such that ${ \phi_\mathfrak{g}}\circ J$ is invariant with respect to the hom-Levi-Civita product, i.e., $L_u\circ { \phi_\mathfrak{g}}\circ J={ \phi_\mathfrak{g}}\circ J\circ L_u$ for any $u\in \mathfrak{g}$.
\end{definition}
Note that condition $(u\cdot({ \phi_\mathfrak{g}}\circ J))(v)=({ \phi_\mathfrak{g}}\circ J)(u\cdot v)$ is equivalent with 
\begin{equation}\label{AM1}
({ \phi_\mathfrak{g}}\circ J)(u)\cdot({ \phi_\mathfrak{g}}\circ J)(v)=({ \phi_\mathfrak{g}}\circ J)(({ \phi_\mathfrak{g}}\circ J)(u)\cdot v),
\end{equation}
and
\begin{equation}\label{AM2}
u\cdot v=-({ \phi_\mathfrak{g}}\circ J)(u\cdot({ \phi_\mathfrak{g}}\circ J)(v)).
\end{equation}
\begin{proposition}
All non-abelian K\"{a}hler hom-Lie algebras of dimension 2 are as $(\mathfrak{g}, [\cdot, \cdot], \widehat{\phi}, \widehat{J}_i, <,>_i)$, $i=1,2$, where $\widehat{\phi}$ is given by (\ref{BB1})
and $\widehat{J}_i$ and  $<,>_i$ have the following matrix presentations:
\begin{equation}\label{LEIL4}
\left\{
\begin{array}{cc}
&\hspace{-2cm}[\widehat{J}_{1}]=
\begin{bmatrix}
a&d\\
h&-a
\end{bmatrix},\ \
[<,>_{1}]=
\begin{bmatrix}
-\frac{d}{a}&1\\
1&\frac{h}{a}
\end{bmatrix},\ \
e_1\cdot e_1=a^2e_1+ade_2,\\
&e_1\cdot e_2=ahe_1-a^2e_2,\ \
e_2\cdot e_1=ahe_1-(a^2+1)e_2,\ \ e_2\cdot e_2=h^2e_1-ahe_2,
\end{array}
\right.
\end{equation}%
\begin{equation}
\left\{
\begin{array}{cc}
&[\widehat{J}_2]=
\begin{bmatrix}
0&d\\
\frac{1}{d}&0
\end{bmatrix},\ \
[<,>_2]=
\begin{bmatrix}
<e_1,e_1>_2&0\\
0&\frac{<e_1,e_1>_2}{d^2}
\end{bmatrix},\\
&e_2\cdot e_1=-e_2,\ \
e_1\cdot e_1=e_1\cdot e_2=0,\ \ e_2\cdot e_2=\frac{1}{d^2}e_1,
\end{array}
\right.
\end{equation}%
where $a,d\neq 0$, $a^2+ ad=-1$.
\end{proposition}
\begin{proof}
Applying Proposition \ref{EProp5}, 
 we must to check that one of the structures determined in Proposition \ref{EProp4} is compatible with these products. We consider two cases as follows:

\textbf{Case 1.} $<e_1, e_2>=0$.\\
In this case, we can consider $<,>_2$, because $<,>_1$, is not a pseudo-Riemannian metric, when $<e_1, e_2>=0$. In this case, the product (\ref{LEI4})
 reduces to
\begin{equation}
e_1\cdot e_1=e_1\cdot e_2=0,\ \ \ e_2\cdot e_1=-e_2,\ \ \ e_2\cdot e_2=\frac{1}{d^2}e_1.
\end{equation}
It is easy to see that (\ref{AM2}) is held for $(\mathfrak{g},
[,],  \widehat{\phi}, \widehat{J}_2, <,>_2)$ with the above
product. So this structure is K\"{a}hler.

\textbf{Case 2.} $<e_1, e_2>\neq0$.\\
In this case we consider $\widehat{J}_1$ and $<,>_1$  given by Proposition \ref{EProp4}. Then, (\ref{LEI5})-(\ref{LEI8}) reduce to
\begin{align}
e_1\cdot e_1&=\frac{a^2<e_1, e_1>^2}{d^2\det [<,>]}e_1+\frac{a<e_1, e_1>^2}{d\det [<,>]}e_2,\label{LEI5'}\\
e_1\cdot e_2&=\frac{-h<e_1, e_1>}{d\det [<,>]}e_1+\frac{a<e_1, e_1>}{d\det [<,>]}e_2,\label{LEI6'}\\
e_2\cdot e_1&=\frac{-h<e_1, e_1>}{d\det [<,>]}e_1-\frac{\det [<,>]-\frac{a}{d}<e_1, e_1>}{\det [<,>]}e_2,\label{LEI7'}\\
e_2\cdot e_2&=\frac{-h^2<e_1, e_1>}{ad\det [<,>]}e_1+\frac{h<e_1, e_1>}{d\det [<,>]}e_2.\label{LEI8'}
\end{align}
Direct calculation together $a^2+hd=-1$ give
\[
-\widehat{J}_{1}(e_1\cdot \widehat{J}_{1}(e_1))=\frac{a^2<e_1, e_1>^2}{d^2\det [<,>]}e_1-\frac{<e_1, e_1>}{\det [<,>]}e_2.
\]
So, condition $e_1\cdot e_1=-\widehat{J}_{1}(e_1\cdot \widehat{J}_{1}(e_1))$ gives $<e_1, e_1>=-\frac{d}{a}$. In this case, $\det[<,>_1]=\frac{1}{a^2}$ and so $[<,>_1]$ and the
 Levi-Civita product reduce to (\ref{LEIL4}). Using it we get also
\[
-\widehat{J}_{1}(e_1\cdot \widehat{J}_{1}(e_2))=ahe_1-a^2e_2=e_1\cdot e_2,
\]
\[
-\widehat{J}_{1}(e_2\cdot \widehat{J}_{1}(e_1))=ahe_1-(a^2+1)e_2=e_2\cdot e_1,
\]
\[
-\widehat{J}_{1}(e_2\cdot \widehat{J}_{1}(e_2))=h^2e_1-ahe_2=e_2\cdot e_2.
\]
So (\ref{AM2}) holds.
\end{proof}
From Corollary \ref{FIN}, we deduce the following
\begin{corollary}
There  exists  no non-abelian K\"{a}hler proper hom-Lie
algebra of dimension 2.
\end{corollary}
\begin{example}
We consider the home-Lie algebra $(\mathfrak{g}, [\cdot,\cdot], \phi,\Omega)$ introduced in Example \ref{IMEX}. We define the metric $<,>$ of $\mathfrak{g}$ as follows
\[
\begin{bmatrix}
A&0&0&0\\
0&\frac{a}{b}A&0&0\\
0&0&A&0\\
0&0&0&\frac{a}{b}A
\end{bmatrix}
.
\]
It is easy to see that  $<\phi(e_i),e_j>=0=<e_i,\phi(e_j)>$, for all $i,j=1,2,3,4$, except 
\begin{align*}
<\phi(e_1),e_1>=&-A=<e_1,\phi(e_1)>,\ \ \ \ \ <\phi(e_2),e_2>=\frac{a}{b}A=<e_2,\phi(e_2)>,\\
<\phi(e_3),e_3>=&-A=<e_3,\phi(e_3)>,\ \ \ \ \ <\phi(e_4),e_4>=\frac{a}{b}A=<e_4,\phi(e_4)>.
\end{align*}
Thus (\ref{SL}) is hold and so $(\mathfrak{g}, [\cdot, \cdot], \phi_{\mathfrak{g}}, <,>)$ is a pseudo-Riemannian hom-Lie algebra. If isomorphism $J$ is determined as
\[
J(e_1)=e_3,\ \ \ J(e_2)=e_4,\ \ \ J(e_3)=-e_1, \ \ J(e_4)=- e_2,
\]
then we have 
\[
J^2(e_1)=J^2(e_2)=J^2(e_3)=J^2(e_4)=-Id_\mathfrak{g}.
\]
Moreover using the above equations we deduce 
\[
\phi^2(e_1)=\phi^2(e_2)=\phi^2(e_3)=\phi^2(e_4)=Id_\mathfrak{g},
\]
and
\begin{align*}
(J\circ\phi)e_1=&-e_3=(\phi\circ J)e_1,\ \ \ (J\circ\phi)e_2=e_4=(\phi\circ J)e_2,\\
(J\circ\phi)e_3=&e_1=(\phi\circ J)e_3,\ \ \ (J\circ\phi)e_4=-e_2=(\phi\circ J)e_4.
\end{align*}
Thus $J$ is an almost complex structure on $(g, [\cdot, \cdot], \phi_{\mathfrak{g}})$.  Also it results that 
\begin{align*}
N(e_i,e_j)=0,\ \ \ \ \forall i,j=1,2,3,4,
\end{align*}
i.e., $J$ is a complex structure on $(\mathfrak{g}, [\cdot,\cdot], \phi)$. It is easy to check that $<({ \phi‎‎_\mathfrak{g}}\circ J)(e_i), ({ \phi‎‎_\mathfrak{g}}\circ J)(e_j)>=<e_i, e_j>$, for all $i, j=1,2,3,4$, and so $(\mathfrak{g}, [\cdot, \cdot], \phi_{\mathfrak{g}}, <,>, J)$ is a Hermitian hom-Lie algebra. Now, we study the K\"{a}hlerian property for this hom-Lie algebra. At first we must to obtain the hom-Levi-Civita product for it. If we denote this product with $\cdot$, then we have $e_i\cdot e_j=\sum_{k=1}^{4}A_{ij}^ke_k$, for all $i, j=1, 2, 3, 4$. From Koszul's formula given by (\ref{Koszul}) we get $<e_1\cdot e_1, \phi(e_i)>=0$, for all $i=1, 2, 3, 4$, which give $e_1\cdot e_1=0$. Again (\ref{Koszul}) gives $<e_1\cdot e_2, \phi(e_i)>=0$, for all $i=1, 2,3, 4$, which imply $A_{12}^1=A_{12}^2=A_{12}^3=A_{12}^4=0$. 
Therefore we deduce $e_1\cdot e_2=0$. But $[e_1, e_2]=e_1\cdot e_2-e_2\cdot e_1$ implies $e_2\cdot e_1=ae_3$. In the similar way we can obtain the following 
\begin{align*}
e_3\cdot e_1=&-be_2,\ \ e_2\cdot e_2=ae_4,\ \  e_2\cdot e_4=-ae_2,\ \ \ e_3\cdot e_4=ae_3,\ \ 
  e_3\cdot e_3=be_4,\ \ e_2\cdot e_3=e_3\cdot e_2=-ae_1,
\end{align*}  
and 
\begin{equation*}
e_1\cdot e_3=e_1\cdot e_4=e_4\cdot e_1=e_4\cdot e_2=e_4\cdot e_3=e_4\cdot e_4=0.
\end{equation*} 
Easily we can see that the hom-Levi-Civita product computed in the above satisfies in (\ref{AM1}). Thus $(\mathfrak{g}, [\cdot,\cdot], \phi_{\mathfrak{g}}, J, <,>)$ is a K\"{a}hle hom-Lie algebra.
\end{example} 
\begin{proposition}\label{304}
Let $(\mathfrak{g},[\cdot, \cdot], \phi_{\mathfrak{g}}, <,>, J)$ be a K\"{a}hler hom-Lie algebra. Then 
 $(\mathfrak{g}, [\cdot, \cdot], { \phi_\mathfrak{g}}, \Omega)$ is a symplectic hom-Lie algebra, where 
\begin{equation}\label{*}
\Omega(u,v)=<({ \phi_\mathfrak{g}}\circ J)u,v>.
\end{equation}
\end{proposition}
\begin{proof}
Applying  (\ref{SL}) and (\ref{*}) we obtain
\begin{align}\label{AM300}
&\Omega([u,v],{ \phi‎‎_\mathfrak{g}}(w))+\Omega([v,w],{ \phi‎‎_\mathfrak{g}}(u))+\Omega([w,u],{ \phi‎‎_\mathfrak{g}}(v))\\
&=-<[u,v],({ \phi‎‎_\mathfrak{g}}\circ J) ({ \phi‎‎_\mathfrak{g}} (w))>-<[v,w],({ \phi‎‎_\mathfrak{g}}\circ J) ({ \phi‎‎_\mathfrak{g}} (u))>-<[w,u],({ \phi‎‎_\mathfrak{g}}\circ J) ({ \phi‎‎_\mathfrak{g}} (v))>\nonumber\\
&=-<u\cdot v,({ \phi‎‎_\mathfrak{g}}\circ J) ({ \phi‎‎_\mathfrak{g}} (w))>+<v\cdot u,({ \phi‎‎_\mathfrak{g}}\circ J) ({ \phi‎‎_\mathfrak{g}} (w))>-<v\cdot w,({ \phi‎‎_\mathfrak{g}}\circ J) ({ \phi‎‎_\mathfrak{g}} (u))>\nonumber\\
&\ \ \ \ +<w\cdot v,({ \phi‎‎_\mathfrak{g}}\circ J) ({ \phi‎‎_\mathfrak{g}} (u))>-<w\cdot u,({ \phi‎‎_\mathfrak{g}}\circ J) ({ \phi‎‎_\mathfrak{g}} (v))>+<u\cdot w,({ \phi‎‎_\mathfrak{g}}\circ J) ({ \phi‎‎_\mathfrak{g}}( v))>,\nonumber
\end{align}
for any $u,v,w\in\mathfrak{g}$. But using (\ref{L2}) and (\ref{AM2})  we conclude 
\begin{align*}\label{AM301}
<u\cdot v, ({ \phi‎‎_\mathfrak{g}}\circ J)({ \phi‎‎_\mathfrak{g}} (w))>=&-<({ \phi‎‎_\mathfrak{g}}\circ J)(u\cdot ({ \phi‎‎_\mathfrak{g}}\circ J)(v)), ({ \phi‎‎_\mathfrak{g}}\circ J)({ \phi‎‎_\mathfrak{g}} (w))>\\
&=-<u\cdot ({ \phi‎‎_\mathfrak{g}}\circ J)v, ({ \phi‎‎_\mathfrak{g}} w)>=<u\cdot w,{ \phi‎‎_\mathfrak{g}} ({ \phi‎‎_\mathfrak{g}}\circ J)(v)>.\nonumber
\end{align*}
Setting the above equation in (\ref{AM300}) we get
\begin{align*}
&\Omega([u,v],{ \phi‎‎_\mathfrak{g}}(w))+\Omega([v,w],{ \phi‎‎_\mathfrak{g}}(u))+\Omega([w,u],{ \phi‎‎_\mathfrak{g}}(v))\\
&=-<u\cdot w, ({ \phi‎‎_\mathfrak{g}}^2\circ J) v>+<v\cdot w, ({ \phi‎‎_\mathfrak{g}}^2\circ J)u>-<v\cdot w,J u>\\
&+<w\cdot u, ({ \phi‎‎_\mathfrak{g}}^2\circ J)v>-<w\cdot u,J v>+<u\cdot w,J v>=0.
\end{align*}
Moreover (\ref{*}) implies
\[
\Omega({ \phi‎‎_\mathfrak{g}}(u),{ \phi‎‎_\mathfrak{g}}(v))=<({ \phi‎‎_\mathfrak{g}}\circ J){ \phi‎‎_\mathfrak{g}}(u), { \phi‎‎_\mathfrak{g}}( v)>=<J(u),{ \phi‎‎_\mathfrak{g}}(v)>=<({ \phi‎‎_\mathfrak{g}}\circ J)u,v>=\Omega(u,v).
\]
\end{proof}
According to the above theorem, a K\"{a}hler hom-Lie algebra has two products, the hom-Levi-Civita product and the hom-left
symmetric product $\bold{ a}$ associated with $(\mathfrak{g}, [\cdot, \cdot], \phi_{\mathfrak{g}}, \Omega)$.
Therefore we have the following:
\begin{corollary}
If $(\mathfrak{g}, [\cdot, \cdot], { \phi_\mathfrak{g}}, <,>, J, \Omega)$ is a K\"{a}hler hom-Lie algebra, then $(\mathfrak{g}, \bold{a},\phi_{\mathfrak{g}})$ is a hom-left symmetric algebra. Also, 
 $(\mathfrak{g}\oplus\mathfrak{g}^*, \bold{a},\Phi)$ is a hom-left symmetric algebra, where $\mathfrak{g}^*$ is the dual space of $\mathfrak{g}$ and $\bold{a}$ and $\Phi$ are given by (\ref{F2}). Moreover,  $(\mathfrak{g}\oplus\mathfrak{g}^*, \Omega_{_{\mathfrak{g}\oplus\mathfrak{g}^*}})$ is a phase space of  hom-Lie algebra $\mathfrak{g}$ where $\Omega_{_{\mathfrak{g}\oplus\mathfrak{g}^*}}$ is given by (\ref{AM505}).
\end{corollary}
From Proposition \ref{Prop4.11} we deduce the following:
\begin{corollary}
Let $(\mathfrak{g}, [\cdot, \cdot], { \phi_\mathfrak{g}}, <,>, J, \Omega)$ be a K\"{a}hler hom-Lie algebra. Then there is a complex structure $\mathcal{J}$ on $\mathfrak{g}\oplus\mathfrak{g}^*$ given by (\ref{F14}).
\end{corollary}


\end{document}